\documentclass[12pt,a4paper]{amsart}
\usepackage{amssymb,cite}
\usepackage{graphicx}
\usepackage{dcpic,pictex}

\theoremstyle{plain}
\newtheorem{theorem}{Theorem}[section]

\newtheorem{lemma}[theorem]{Lemma}

\newtheorem{conjecture}{Conjecture}

\theoremstyle{definition}

\theoremstyle{remark}

\numberwithin{equation}{section}

\numberwithin{table}{section}

\numberwithin{figure}{section}

\begin{document}

\title[The $\mu$-permanent revisited]{The $\mu$-permanent revisited}

%    Information for first author
\author{Carlos M. da Fonseca}
\address{University of Primorska, FAMNIT, Glagolja\v ska 8, 6000 Koper, Slovenia}
\email{carlos.dafonseca@famnit.upr.si}

\thanks{This manuscript is largely based on the talk ``On a generalization of the determinant of a matrix" that the author gave in Washington \& Lee University, Lexington, VA, USA, on July 19, 2012. It was submitted to Linear and Multilinear Algebra on September 20, 2015. A reduced form will be published as a Letter to the Editor in the same journal.}

%    General info

\date{\today}

%\date{October 8, 2015}


\begin{abstract}
Let $A=(a_{ij})$ be an $n$-by-$n$ matrix. For any real number $\mu$,
we define the polynomial
$$P_\mu(A)=\sum_{\sigma\in S_n}
a_{1\sigma(1)}\cdots a_{n\sigma(n)}\,\mu^{\ell(\sigma)}\; ,$$ as the
$\mu$-permanent of $A$, where $\ell(\sigma)$ is the number of
inversions of the permutation $\sigma$ in the symmetric group $S_n$.
In this note, we review several less known results of the
$\mu$-permanent, recalling some of its interesting properties. Some
determinantal conjectures are considered and extended to that
polynomial. A correction to a previous note is presented as well.
\end{abstract}

\maketitle

\section{Introduction}

Given an $n\times n$ matrix $A=(a_{ij})$ and a real number $\mu$, we
define the $\mu$-permanent of $A$ as the polynomial
\begin{equation}\label{mup}
P_\mu(A)=\sum_{\sigma\in S_n} \left(\prod_{i=1}^{n}
a_{i\sigma(i)}\right)\mu^{\ell(\sigma)}\; ,
\end{equation}
where $\ell(\sigma)$ is the number of inversions of the permutation
$\sigma$ in the symmetric group $S_n$ of degree $n$, i.e., the
number of interchanges of consecutive elements necessary to arrange
$\sigma$ in its natural order \cite[p.1]{Muir}, i.e.,
$$\ell(\sigma)=\#\{(i,j)\in \{1,\ldots,n\}^2\, |\, i<j \mbox{ and }
\sigma(i)>\sigma(j)\}\, .$$ For example,
$$
  P_\mu
\left(
  \begin{array}{cc}
    a_{11} & a_{12}   \\
    a_{21} & a_{22}
  \end{array}
\right)
   = a_{11}a_{22}+a_{12}a_{21}\mu
$$
and
\begin{eqnarray*}
  P_\mu
\left(
  \begin{array}{ccc}
    a_{11} & a_{12} & a_{13} \\
    a_{21} & a_{22} & a_{23} \\
    a_{31} & a_{32} & a_{33} \\
  \end{array}
\right)
   &=& a_{11}a_{22}a_{33}+a_{11}a_{23}a_{32}\,\mu+a_{12}a_{21}a_{33}\,\mu+ \\
   && +\,
   a_{12}a_{23}a_{31}\, \mu^2+a_{13}a_{21}a_{32}\,\mu^2+a_{13}a_{22}a_{31}\,\mu^3\,
   .
\end{eqnarray*}

The $\mu$-permanent of a square matrix seems fairly trivial to
manipulate, but it always leads into lengthy and tedious
calculations, and it is often notoriously hard to compute
\cite[p.190]{BR}.

Among the linear algebra community, this notion was first introduced
in $1991$ by Ravindra Bapat \cite{B} as a possible interpolation of
the determinant ($\mu=-1$) and the permanent ($\mu=1$). Bapat
\emph{et al.} also called it later as the $q$-permanent of $A$ and
wrote per$_q(A)$ \cite{BL,BR,L2}). Actually, this has become the
most common designation for it. However, as we will see, by the same
time, this notion emerged independently in other fields of
mathematics: Grassmann algebras and quantum groups. In any case, it
should not be confused with the $\alpha$-determinant defined by D.
Vere-Jones \cite{VJ}, where $\ell(\sigma)$ is $n$ minus the number
of cycles in $\sigma\in S_n$.

In this note, we clarify this notion presenting some less-known
basic results. We correct a typo and some consequences in our
previous paper \cite{dF2005}. We also discuss some conjectures on
determinants and permanents for the $\mu$-permanent and prove them
for particular families of matrices. A more recent conjecture for
the $\mu$-permanent of a tridiagonal matrix is also considered.

\section{The \textbf{q}-determinant}

In June, $1989$, Yang submited a note \cite{Y1991} where he defined
the $q$-determinant of a matrix $A$, $\det_q(A)$, in the same way as
in \eqref{mup}, where the entries are in a certain commutative ring.
No particular previous motivation was given. The main aim was to
extend the analysis of the determinant to $q$-Grassmann algebras.
Briefly, a $q$-Grassmann algebra is the associative
$\mathbb{K}[q]-$algebra, where $\mathbb{K}$ is a field of
characteristic $0$, generated by $x_1,\ldots,x_n$, satisfying the
relations $x_i^2=0$ and $x_ix_j = qx_jx_i$, if $i < j$. A Grassmann
algebra is also known as exterior algebra.

Coincidentally -  or not - exactly at the same time,  Noumi, Yamada,
and Mimachi, led by representations of the quantum groups, announced
in \cite{NYM1989} the same definition where $q$ is replaced by $-q$.
Here the quantum determinant is defined over a $\mathbb{C}$-algebra
where the canonical generators satisfied certain relations,
containing the relations studied by Yang. Later on, the same authors
found several useful properties in \cite{NYM1993}, namely the
Laplace expansion. To the best of our knowledge, this is the first
time the concept of $\mu$-permanent was defined. Its motivation is
clearly independent from the one of Bapat and Lal (cf.
\cite{B,BL,BR,L,L2}).

Probably the main result one can find in \cite{Y1991} is Theorem
1(4), where the Laplace expansion through the first/last row/column
is stated. Namely,
\begin{eqnarray*}
% \nonumber to remove numbering (before each equation)
   {\rm det}_{q}(A)    &=&  \sum_{j=1}^nq^{j-1} a_{1j}\,{\rm det}_q(A_{1j})=\sum_{j=n}^1q^{n-j}a_{nj}\, {\rm det}_q(A_{1j})\\
                &=& \sum_{i=1}^nq^{i-1}a_{i1}\,{\rm det}_q(A_{i1})=\sum_{i=n}^1q^{n-i}a_{in}\, {\rm
                det}_q(A_{in})\, ,
\end{eqnarray*}
where $A_{ij}$ is the $(n-1)\times(n-1)$ submatrix obtained from
deleting the $i$th row and the $j$th column of $A$. This observation
is very important to avoid confusion in the future. Moreover, the
original definition of $A_{ij}$ as the ``$(i,j)$-minor of $A$" in
\cite[Theorem 1(4)]{Y1991} is not very accurate.

Shortly after, Tagawa \cite{T1991} considered the same concept for
commutative rings. Basically, the definition was the same as
\eqref{mup}. To be more precise, we have $\det_{-\mu}(A)=P_\mu(A)$.
This author \cite{T1993} designated exactly $\det_{-q}(A)$ by
$q$-permanent denoting it by ${\rm per}_q(A)$. Furthermore, Tagawa
extended even more the concept defining the multivariable quantum
determinant. For a sequence $\textbf{q}= (q_{1},
q_{2},\ldots,q_{n})$, the multivariable quantum determinant, or
simply \textbf{q}-determinant, of $A$ is defined by
$$
{\rm det}_\textbf{q}=\sum_{\sigma\in S_n} \prod_{i=1}^{n}
a_{i\sigma(i)}(-q_i)^{\ell_i(\sigma)}\; ,
$$
where $\ell_i(\sigma)$ is the number of inversions of $\sigma$ at
$i$. For example, we have
$$
{\rm det}_\textbf{q} \left(
  \begin{array}{cc}
    a_{11} & a_{12}  \\
    a_{21} & a_{22}
  \end{array}
\right)= a_{11}a_{22}-a_{12}a_{21}q_1
$$
and
\begin{eqnarray*}
{\rm det}_\textbf{q} \left(
  \begin{array}{ccc}
    a_{11} & a_{12} & a_{13} \\
    a_{21} & a_{22} & a_{23} \\
    a_{31} & a_{32} & a_{33} \\
  \end{array}
\right) & = &
a_{11}a_{22}a_{33}-a_{11}a_{23}a_{32}q_2-a_{12}a_{21}a_{33}q_1 \\
   & & \; + \, a_{12}a_{23}a_{31}q_1q_2+a_{13}a_{21}a_{32}q_1^2-a_{13}a_{22}a_{31}q_1^2q_2.
\end{eqnarray*}

Clearly, from the definition, the determinant and the permanent as
well as the $q$-determinant are particular cases of the
\textbf{q}-determinant. For the \textbf{q}-determinant, the
multilinearities with respect to the rows and the columns are valid
as in the case of the ordinary determinant. This allows us to state
the same for the $q$-determinant or the $\mu$-permanent. However, in
general, $\det_\textbf{q}(A^T)\neq\det_\textbf{q}(A)$, but
$\det_{q}(A^T)=\det_{q}(A)$ is always true. These observations can
be found in \cite{T1993,T1991,Y1991}.

We find in \cite{T1991} the first properties for the $q$-determinant
as particular cases of the corresponding relation involving the
\textbf{q}-determinant and later extended in contemporary papers
\cite{NYM1993,T1993}. We will translate next some of them to the
$\mu$-permanent. If $A=(a_{ij})$ is an $n$-by-$n$ matrix, we define
the $(i,j)$-$q$-complementary matrix $A_{ij}(q)$ of $A$ by
$$
A_{ij}(q)=
\left(
  \begin{array}{ccc|ccc}
    a_{1,1} & \cdots & a_{1,j-1} & qa_{1,j+1} & \cdots & qa_{1,n} \\
    \vdots & \vdots & \vdots & \vdots & \vdots & \vdots \\
     a_{i-1,1} & \cdots & a_{i-1,j-1} & qa_{i-1,j+1} & \cdots & qa_{i-1,n}
     \\ \hline
     qa_{i+1,1} & \cdots & qa_{i+1,j-1} & a_{i+1,j+1} & \cdots & a_{i+1,n}  \\
     \vdots & \vdots & \vdots & \vdots & \vdots & \vdots \\
     qa_{n,1} & \cdots & qa_{n,j-1} & a_{n,j+1} & \cdots & a_{n,n}
  \end{array}
\right)\, .
$$

Then, from \cite[Lemma 2.3]{T1991}, we have
\begin{equation} \label{delete}
P_\mu(A_{ij}(\mu))= P_\mu \left(
  \begin{array}{ccc|c|ccc}
    a_{1,1} & \cdots & a_{1,j-1} & 0 & a_{1,j+1} & \cdots & a_{1,n} \\
    \vdots & \vdots & \vdots & \vdots & \vdots & \vdots & \vdots \\
     a_{i-1,1} & \cdots & a_{i-1,j-1} & 0 & a_{i-1,j+1} & \cdots & a_{i-1,n}
     \\ \hline
    0 & \cdots & 0 & 1 & 0 & \cdots & 0
     \\ \hline
     a_{i+1,1} & \cdots & a_{i+1,j-1} & 0 & a_{i+1,j+1} & \cdots & a_{i+1,n}  \\
     \vdots & \vdots & \vdots & \vdots & \vdots & \vdots & \vdots \\
     a_{n,1} & \cdots & a_{n,j-1} & 0 & a_{n,j+1} & \cdots & a_{n,n}
  \end{array}
\right) .
\end{equation}
Therefore we have the Laplace expansion (cf.
\cite[Proposition 1.2.2]{T1993} or \cite[Proposition 2.2]{T1991})
\begin{equation}\label{laplace}
    P_\mu(A)=\sum_{j=1}^n a_{ij}P_\mu(A_{ij}(\mu))=\sum_{i=1}^n
    a_{ij}P_\mu(A_{ij}(\mu))\, .
\end{equation}

Observe that \eqref{laplace} coincides with the corollary of
Proposition 1.1 in \cite{NYM1993}, due to the fundamental relations
(1.1) of the generators of the $\mathbb{C}$-algebra in that paper.

For another approach to quantum determinant the reader is referred
to \cite{HH1992}.

\section{Corrections of previous results}

It is clear that, under similarity, the $\mu$-permanent does not
keep the same value, in general,  i.e., the polynomial $P_\mu(A)$ is
not necessarily the same as $P_\mu(BAB^{-1})$, for $B$ nonsingular.
In particular, for permutation similarity, this means that
interchanging rows and columns of the same indices leads to possible
different $\mu$-permanents. Indeed, if $P(\tau)$ is the permutation
matrix $(\delta_{i,\tau(i)})$, then
$$
  P_\mu(P(\tau)^{-1}AP(\tau)) = \sum_{\sigma\in S_n} \prod_{i=1}^{n}
a_{i\sigma(i)}\,\mu^{\ell(\tau\sigma\tau^{-1})}\, .
$$

In the note \cite{dF2005}, we first analyzed the $\mu$-permanent of
a matrix in terms of the underlying graph of the matrix, with
special focus on the (symmetric) matrices whose graph is a tree.
Interchanging rows and columns does not change the the underlying
graph (or digraph) of the matrices involved. However the labeling of
the vertices changes. Saying that, some corrections are needed.

First we need to clarify a notion in \cite{dF2005} which is
misleading. When we have a digraph $D$, in the sense of
\eqref{delete}, $D\backslash i$ is obtained by deleting all arcs
including $i$, but not exactly the vertex $i$ as it was stated. This
is implicit, for example, in the discussion of the example in
Section 5. The existence of isolated vertices do not interfere in
any result. So, for example, in the definition of a tree, we have a
connected acyclic graph with eventually isolated vertices.

Next, we need an additional condition on the labeling of all graphs:
given two disjoint edges $ij$ and $k\ell$ (say $i<j$, $k<\ell$, and
$i<k$), then one of the conditions is satisfied:
\begin{itemize}
  \item[(i)] $i<j<k<\ell$
  \item[(ii)] $i<k<\ell<j$.
\end{itemize}

Clearly, for any tree, such labeling is always possible. For
example, the results are valid for any matrices whose graph is

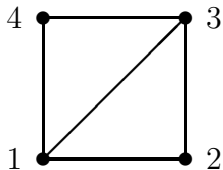
\begin{center}
\setlength{\unitlength}{0.75mm} \thicklines
\begin{picture}(35,30)

\put(5,0){\circle*{2.5}} \put(30,0){\circle*{2.5}}
\put(5,25){\circle*{2.5}} \put(30,25){\circle*{2.5}}

\put(5,0){\line(1,1){25}} \put(5,0){\line(1,0){25}}
\put(5,0){\line(0,1){25}}

\put(5,25){\line(1,0){25}} \put(30,0){\line(0,1){25}}

\put(0,0){\makebox(0,0){$1$}} \put(35,0){\makebox(0,0){$2$}}
\put(0,25){\makebox(0,0){$4$}} \put(35,25){\makebox(0,0){$3$}}

\end{picture}
\end{center}

\hspace{1em}

For paths, we have for example

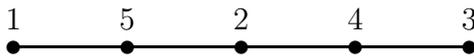
\begin{center}
\setlength{\unitlength}{0.75mm} \thicklines
\begin{picture}(80,15)

\put(0,0){\circle*{2.5}} \put(20,0){\circle*{2.5}}
\put(40,0){\circle*{2.5}} \put(60,0){\circle*{2.5}}
\put(80,0){\circle*{2.5}}

\put(0,0){\line(1,0){80}}

\put(0,5){\makebox(0,0){$1$}} \put(20,5){\makebox(0,0){$5$}}
\put(40,5){\makebox(0,0){$2$}} \put(60,5){\makebox(0,0){$4$}}
\put(80,5){\makebox(0,0){$3$}}

\end{picture}
\end{center}

\hspace{1em}

or

\begin{center}
\setlength{\unitlength}{0.75mm} \thicklines
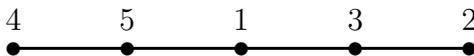
\begin{picture}(80,15)

\put(0,0){\circle*{2.5}} \put(20,0){\circle*{2.5}}
\put(40,0){\circle*{2.5}} \put(60,0){\circle*{2.5}}
\put(80,0){\circle*{2.5}}

\put(0,0){\line(1,0){80}}

\put(0,5){\makebox(0,0){$4$}} \put(20,5){\makebox(0,0){$5$}}
\put(40,5){\makebox(0,0){$1$}} \put(60,5){\makebox(0,0){$3$}}
\put(80,5){\makebox(0,0){$2$}}

\end{picture}
\end{center}

\hspace{1em}

or, for more general trees,

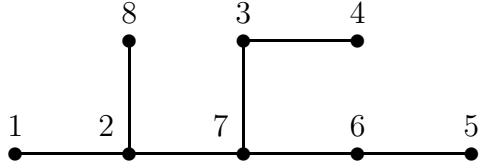
\begin{center}
\setlength{\unitlength}{0.75mm} \thicklines
\begin{picture}(80,30)

\put(0,0){\circle*{2.5}} \put(20,0){\circle*{2.5}}
\put(20,20){\circle*{2.5}} \put(40,0){\circle*{2.5}}
\put(40,20){\circle*{2.5}} \put(60,0){\circle*{2.5}}
\put(60,20){\circle*{2.5}} \put(80,0){\circle*{2.5}}

\put(0,0){\line(1,0){80}} \put(20,0){\line(0,1){20}}
\put(40,0){\line(0,1){20}} \put(40,20){\line(1,0){20}}

\put(0,5){\makebox(0,0){$1$}} \put(16,5){\makebox(0,0){$2$}}
\put(20,25){\makebox(0,0){$8$}} \put(36,5){\makebox(0,0){$7$}}
\put(60,5){\makebox(0,0){$6$}} \put(80,5){\makebox(0,0){$5$}}
\put(40,25){\makebox(0,0){$3$}} \put(60,25){\makebox(0,0){$4$}}

\end{picture}
\end{center}

\hspace{1em}

It is important now to clarify that Corollary 3.3, Lemma 4.1, and
Theorem 4.3 in \cite{dF2005} are valid whenever the edges of the
underlying graph of the matrix satisfy either (i) or (ii).

Finally, a large number of computational experiments for several
families of graphs lead us to conjecture (cf. Conjecture
\ref{conjx}) that $P_\mu(A)$ is strictly increasing for $\mu\in
[-1,\infty)$, assuming that $A$ is any positive definite, extending
\cite[Theorem 4.3]{dF2005}.

\section{A conjecture over the positive definiteness}

We have seen that the $\mu$-permanent is a  parametric
generalization of both the determinant and the permanent, making
$\mu=-1$ and $\mu=1$, respectively. Note also that
$P_0(A)=a_{11}\cdots a_{nn}$.

The Schur power matrix of $A$, denoted by $\Pi(A)$ \cite{Merris,S2}, is the $n!\times n!$ matrix whose rows and
columns are indexed by $S_n$, where the $(\sigma,\tau)$-entry is
$$\prod_{i=1}^n a_{\sigma(i)\tau(i)}\, ;$$ $\Gamma_\mu$ is a matrix of same
type with the $(\sigma,\tau)$-entry defined as
$$\mu^{\ell(\tau\sigma^{-1})}.$$ For more results and several open problems on $\Gamma_\mu$ the reader is referred to
\cite{C,Stanciu,Z}.
If we set
$$\Pi_\mu(A)=\Pi(A)\circ\Gamma_\mu\; ,$$ where $\circ$ denotes the
Hadamard product, then we have
\begin{equation}\label{Pmu}
    P_\mu(A)=\frac 1{n!}\, \langle\Pi_\mu(A)
    \mathbf{1},\mathbf{1}\rangle ,
\end{equation}
where $\mathbf{1}$ denotes the column vector of all ones.

Bo\.zejko and Speicher \cite{BS} proved that $\Gamma_\mu$ is
positive semidefinite, for $\mu\in[-1,1]$, and it is known that if
$A$ is positive semidefinite matrix, then $\Pi(A)$ is also positive
semidefinite \cite{BL}. Therefore, from (\ref{Pmu}), we
may state:

\begin{lemma}\cite{B,BL} \label{lem}
For any Hermitian positive semidefinite matrix $A$,
$$P_\mu(A)\geqslant 0\; , \quad \mbox{ if }\; \mu\in[-1,1]\;.$$
\end{lemma}

Remark that, although the assumption of the complex (semi)definite matrix being Hermitian is unnecessary, we will include it, since we often consider real matrices, and in that case the symmetry assumption should be incorporated.

Bapat and Lal also established several conjectures in \cite{B,BL}. One of them is the following:

\begin{conjecture}\cite{B,BL} \label{conj1}
Given an $n\times n$ Hermitian positive definite nondiagonal matrix
$A$, $P_\mu(A)$ is a strictly increasing function of $\mu\in[-1,1]$.
\end{conjecture}

This conjecture was motivated by the classical Hadamard inequality
and the permanental analogue proved by Marcus \cite{M}. It has been
proved for $n\leqslant 3$ in \cite{B}, for tridiagonal positive
definite matrices in \cite{L}, and in \cite{dF2005} for symmetric
positive definite matrices under the graph labeling discussed in the
previous section.

More recently, Conjecture \ref{conj1} was extended \cite{dF2}.

\begin{conjecture}\cite{dF2} \label{conjx}
For a given matrix $A>0$, there exists $\epsilon< -1$ such that
$P_\mu(A)$ is a strictly increasing function of $\mu\in
(\epsilon,+\infty)$.
\end{conjecture}

Using orthogonal polynomials and chain sequences the author proved Conjecture \ref{conjx} for tridiagonal matrices \cite{dF2}. In this section we prove this conjecture for other classes of matrices. We start with the case of a $3\times 3$ positive definite matrix
\begin{equation}\label{3by3}
A=
\left(
  \begin{array}{ccc}
    a_{11} & a_{12} & a_{13} \\
    a_{12} & a_{22} & a_{23} \\
    a_{13} & a_{23} & a_{33} \\
  \end{array}
\right)\, .
\end{equation}

From \eqref{mup}, we have
$$
  P_\mu(A) = a_{11}a_{22}a_{33}+a_{11}a_{23}^2\,\mu+a_{12}^2a_{33}\,\mu+
  2\, a_{12}a_{13} a_{23}\, \mu^2+a_{13}^2a_{22}\,\mu^3\, .
$$
Therefore
\begin{eqnarray*}
   \frac{d}{d\mu}\, P_\mu(A) &=&a_{11}a_{23}^2+a_{12}^2a_{33}+
  4\, a_{12}a_{13} a_{23}\, \mu+ 3\, a_{13}^2a_{22}\,\mu^2 \\
   &=& \left(
  \begin{array}{cc}
    a_{23} & a_{13}\,\mu
  \end{array}
\right) \left(
  \begin{array}{cc}
    a_{11} & a_{12}  \\
    a_{12} & a_{22} \\
  \end{array}
\right) \left(
  \begin{array}{c}
    a_{23} \\ a_{13}\,\mu
  \end{array}
\right)  + \\
   && \quad + \, \left(
  \begin{array}{cc}
    a_{13}\,\mu & a_{12}
  \end{array}
\right) \left(
  \begin{array}{cc}
    a_{22} & a_{23}  \\
    a_{23} & a_{33} \\
  \end{array}
\right) \left(
  \begin{array}{c}
    a_{13}\, \mu \\ a_{12}
  \end{array}
\right) +\\
&& \qquad + \,
a_{13}^2a_{22}\,\mu^2
\end{eqnarray*}
is always positive, which implies $P_\mu(A)$ is a strictly increasing function. Since $A$ is positive definite, $P_{-1}(A)>0$ and,
consequently, the only zero of $P_\mu(A)$ is less than $-1$. Moreover
$$P_\mu(A)> 0\; , \quad \mbox{ for any }\; \mu\in(-1,\infty)\;.$$

A similar conclusion can be reached, for example, for matrices whose graph is a star, with the central vertex labeled with $1$, i.e.,
$$A=
\left(
  \begin{array}{ccccc}
    a_{11} & a_{12} & a_{13} & \cdots & a_{1n}\\
    a_{12} & a_{22} & 0 & \cdots & 0\\
    a_{13} & 0 & a_{33} & \ddots & \vdots\\
    \vdots & \vdots & \ddots & \ddots & 0 \\
     a_{1n} & 0 & \cdots & 0 & a_{nn} \\
  \end{array}
\right)\, .
$$
From \cite[Corollary 3.3]{dF2005}
$$
  P_\mu(A) = a_{11}a_{22}\cdots a_{nn}+\sum_{k=2}^n  a_{1k}^2\left(\prod_{1\neq i\neq k}^n a_{ii}\right) \mu^{2k-3}\, .
$$
Hence
$$ \frac{d}{d\mu}\, P_\mu(A) = \sum_{k=2}^n (2k-3)a_{1k}^2 \left(\prod_{1\neq i\neq k}^n a_{ii}\right)\mu^{2k-4}$$
which is always positive, for any $\mu$.

\section{Some inequalities}

The next open problem is based on the analogue of Lieb's inequality
(cf. \cite{Lieb}) to $\mu$-permanent.

\begin{conjecture}\cite{BL} \label{conj2}
Given an $n\times n$ Hermitian positive semidefinite matrix $A$, let
$S$ be a nonempty subset of $\{1,\ldots,n\}$. Then for $\mu\in [0,1]$,
\begin{equation}\label{liebs}
    P_\mu(A)\geqslant \sum_{{\sigma\in S_n}\atop {\sigma (S)=S}} \prod_{i=1}^{n}
a_{i\sigma(i)}\,\mu^{\ell(\sigma)}\; .
\end{equation}
\end{conjecture}

Conjecture \ref{conj2} can be easily verified for matrices whose
graph is a tree under the labeling previously established. In fact,
since $A$ is Hermitian positive semidefinite ($n\times n$) and
$\mu\in [0,1]$, the sum (\ref{mup}) is nonnegative and we get
$$\sum_{\sigma\in S_n}\prod_{i=1}^{n}
a_{i\sigma(i)}\,\mu^{\ell(\sigma)}- \sum_{{\sigma\in S_n}\atop
{\scriptscriptstyle\sigma (T)=T}} \prod_{i=1}^{n}
a_{i\sigma(i)}\,\mu^{\ell(\sigma)}=\sum_{{\sigma\in S_n}\atop
{\scriptscriptstyle\sigma (T)\neq T}} \prod_{i=1}^{n}
a_{i\sigma(i)}\,\mu^{\ell(\sigma)}\geqslant 0\; .$$

Lal \cite{L} proved the simplest version of Conjecture \ref{conj2}
which is in fact a generalization of Fischer's inequality for
determinants:

\begin{theorem}\cite{L} \label{theo1}
Let $A$ be an $n\times n$ Hermitian positive semidefinite matrix
which is partitioned as
$$A=
\left(
  \begin{array}{cc}
    A_{11} & A_{12} \\
    A^*_{12} & A_{22} \\
  \end{array}
\right)
$$ where $A_{11},A_{22}$ are Hermitian. If
$\mu\in [0,1]$, then
\begin{equation}\label{liebs2}
   P_\mu(A)\geqslant P_\mu(A_{11})P_\mu(A_{22})\; .
\end{equation}
\end{theorem}

Conjecture \ref{conj2} is also true for $n=3$.  We only have to
check the case when $S=\{ 2\}$. The other cases are applications of
Theorem \ref{theo1}. So, we have
$$a_{11}|a_{23}|^2+|a_{12}|^2a_{33}+
   \bar{a}_{12}a_{13}\bar{a}_{23}\mu+a_{12}\bar{a}_{13}a_{23}\mu\geqslant
   0.$$
In fact,
$$
\left(
  \begin{array}{cc}
   \bar{a}_{23} & a_{12}
  \end{array}
\right)
\left(
  \begin{array}{cc}
    a_{11} & a_{13}\mu \\
    \bar{a}_{13}\mu &  a_{33} \\
  \end{array}
\right)
\left(
  \begin{array}{c}
    a_{23}  \\
    \bar{a}_{12}
  \end{array}
\right)\geqslant 0.
$$

Bapat and Lal also extended the longstanding Soules' conjecture
\cite{Merris,S,S2} to a $\mu$-per\-ma\-nental form.

\begin{conjecture}\cite{BL} \label{conj3}
Given an $n\times n$ Hermitian positive semidefinite matrix $A$, let
$\mu\in [0,1]$. Then the largest eigenvalue of $\Pi_\mu(A)$ is
$P_\mu(A)$.
\end{conjecture}

Conjecture \ref{conj3} has been proved for $n\leqslant 3$ in
\cite{L}. For higher orders no progress has been made so far.
However, for acyclic matrices, Conjecture \ref{conj3} can be easily
shown. In fact, since $\Pi_\mu(A)$ is a Hermitian positive
semidefinite $n\times n$ matrix, all eigenvalues are real and
nonnegative. The eigenvectors associated to different eigenvalues
are therefore orthogonal. Note also that, since the graph of $A$ is
a tree, in each column (and in each row) of $\Pi_\mu(A)$, we have
all terms of $P_\mu(A)$. Moreover, each term of the sum (\ref{mup})
appears in each column (and in each row) of $\Pi_\mu(A)$ only once.

From (\ref{Pmu}), $\mathbf{1}$ is an eigenvector of $\Pi_\mu(A)$
associated to the eigenvalue $P_\mu(A)$. Supposing that $\mathbf{u}$
is an eigenvector of $\Pi_\mu(A)$ associated to the eigenvalue
$\lambda(\neq P_\mu(A))$, then $$u_1+\cdots+u_n=0\; .$$ Let us
assume that $u_k$ is largest (positive) coordinate of $\mathbf{u}$.
Since
$$\Pi_\mu(A)\,\mathbf{u}=\lambda\,\mathbf{u}\; ,$$ we have
$$\sum_{j=1}^n\Pi_\mu(A)_{kj}\,u_j=\lambda\,u_k\; ,$$ and, finally, $$
\lambda\,u_k=\sum_{j=1}^n\Pi_\mu(A)_{kj}\,u_j\leqslant
\sum_{j=1}^n\Pi_\mu(A)_{kj}\,u_k=P_\mu(A)\, u_k\; ,$$ because all
entries of $\Pi_\mu(A)$ are nonnegative. Hence $$\lambda\leqslant
P_\mu(A)\; .$$ Thus, we may state:

\begin{theorem}
Given an $n\times n$ Hermitian positive semidefinite matrix $A$
whose graph is a tree, let $\mu\in [0,1]$. Then the largest eigenvalue
of $\Pi_\mu(A)$ is $P_\mu(A)$.
\end{theorem}

Again, the previous theorem is established for trees with the
labeling satisfying the discussed conditions.

\section{Final remarks}

In \cite{BL}, Bapat and Lal started the analysis of a
$\mu$-permanental analogue of Gram's equality. Later on, Lal
\cite{L2} considered the Gram's inequality for the Schur power
matrix, inequalities using induced matrices, and inequalities of
Schwarz type providing generalizations to the $\mu$-permanent of
some results by Ando \cite{A}, and by Marcus and Minc \cite{MM}.
Note that some inequalities of Minkowski type presented by Ando
\cite{A} can also be extended to the $\mu$-permanent.

\end{document}